\documentclass[12pt]{article}

\usepackage{amsmath,amsthm,amssymb}

\usepackage[cp1251]{inputenc}
\usepackage[english]{babel}

\usepackage{authblk}

\usepackage{cite}
\usepackage[colorlinks,linkcolor=blue,citecolor=blue,filecolor=blue,urlcolor=blue]{hyperref}
\usepackage[ruled,vlined,linesnumbered]{algorithm2e}
\usepackage{xcolor}

\usepackage{fancybox}
\usepackage{graphicx}

\theoremstyle{plain}

\theoremstyle{definition}

\theoremstyle{remark}
\newtheorem{remark}{Remark}

\def\Q{\widehat{Q}}

\begin{document}

\title{Discrete orthogonal polynomials \\ as a tool for detection of small anomalies of time series: \\
a case study of GPS final orbits}

\author[1]{S.~P.~Tsarev \thanks{E-mail: \href{mailto:sptsarev@mail.ru}{sptsarev@mail.ru}}}
\author[1]{A.~A.~Kytmanov \thanks{E-mail: \href{mailto:aakytm@gmail.com}{aakytm@gmail.com}}}

\affil[1]{School of Space and Information Technology, Siberian Federal University, 79~Svobodny~pr., 660041 Krasnoyarsk, Russia}

\date{}

\maketitle

\begin{abstract}
In this paper, we present the following results:
\begin{itemize}
    \item We show that the classical discrete orthogonal univariate polynomials (namely, Hahn polynomials on an equidistant lattice with unit weights) of sufficiently high degrees have extremely small values near the endpoints (we call this property as ``rapid decay near the endpoints of the discrete lattice'', cf.~Section~\ref{Fast_decay}).
    \item We demonstrate the importance of the proved results applying polynomial least squares approximation for detection of anomalous values in IGS final orbits for GPS and GLONASS satellites (Section~\ref{sec-lsapprox-IGS}).
    \item We propose a numerically stable method for construction of discrete orthogonal polynomials of high degrees. It allows one to reliably construct Hahn-Chebyshev polynomials using standard accuracy (double precision, 8-byte) on thousands of points, for degrees up to several hundred (Appendix~\ref{App-Hahn-computation}). 
     A Julia implementation of the mentioned algorithms is available at\\
     \url{https://github.com/sptsarev/high-deg-polynomial-fitting}.
\end{itemize}

These results seem to be new; their explanation in the framework of the well-known asymptotic theory of discrete orthogonal polynomials could not be found in the literature \cite{Szego}--\cite{OLB10}.
\end{abstract}

\medskip

\noindent\textbf{Keywords:} discrete orthogonal polynomials; Hahn polynomials; GPS; GLONASS; time series anomalies

\section{Introduction}
\label{sec-intro}


The practical motivation for the theoretical research exposed in this paper was the need to detect small discontinuities in orbits of GPS satellites published by IGS (International GNSS Service, \cite{1}) as ``final GPS and GLONASS
satellite ephe\-me\-ri\-des''. The data files of GPS and GLONASS ephemerides are published by IGS as plain text files organized in a special SP3 format.
These files give the positions of the satellites in the ECEF coordinate systems ("earth-centered, earth-fixed", also known as ECR "earth-centered rotational" cartesian coordinate system) for each day with 15-min time step  and claimed accuracy (RMS)
$\sim$~3~cm. 
The satellite positions are averaged over processing results of several IGS Analysis Centers, according to the nominal IGS standard the published orbit deviates from the calculated average orbit not more than 0.5~mm.
Note that 15 min time steps for satellite positions correspond to $\sim$~3500~km of satellite movement along the orbit, which is
(approximately) circle of radius 25000~km in the inertial non-rotational Earth-centered cartesian coordinate system.
So the problem of detection of small anomalies 
in such fairly sparse time series with large variations seems to be difficult.

As shown below, in fact we can detect the following anomalies in the published GPS orbits (for example, for the period 01.01.2010--01.01.2019):
\begin{itemize}
 \item  ``jumps'' (discontinuities)
 on the day boundaries (at 00:00:00 of every day) with magnitude $\sim $ 1
 cm (this is approximately $10^{-9}$ of the typical values of the GPS coordinates!);
 \item  anomalous big ``jumps'' (up to 100 m)
 on 00:00:00 for some days in the given period;
 \item  ``anomalous singular values'' (usually called ``outliers'')
 inside the day period for very few days; some of them are related to GPS satellite
 maneuvers;
\end{itemize}
Analogous anomalies are found also for the published GLONASS orbits. The detailed account will be presented in a later publication.

The possibility to detect outliers of very small magnitude (approx. $10^{-5}$~km compared to the values of the time series of order $2\cdot 10^{4}$~km) is due to initial smoothness of the trajectories calculated numerically by the IGS Analysis Centers with high precision, and the relatively low rate of the anomalous phenomena described above. In fact, the problem of small (and not always small!) discontinuities on the day boundaries for the published GPS final orbits is well known, a recommended tool for detection and estimation of such discontinuities is numerical modelling of GPS trajectories using numerical integration of their equations of motion with very detailed account of all the forces acting on the satellites, earth orientation parameters and other numerous details of satellite motion modelling \cite{GriffithsRay2009}. In contrast to this complicated anomaly detection method, we propose a very straighforward and methodologically very simple method of detection and estimation of these discontinuities avoiding complete and difficult satellite orbit computations. For this purpose, we use the classical tool of least squares polynomial fitting. The only point one has to take into consideration is the need to find a robust procedure for fitting by polynomials of very high degree (up to 200). The present paper describes the necessary theoretical and computational details of our approach.

The mathematical tool we used in our study is the well-known theory of discrete orthogonal polynomials. In this paper, we establish a number of important new results on asymptotic behavior of the classical Hahn (Chebyshev) discrete orthogonal polynomials (we call this property ``rapid decay near the endpoints'', cf.~Section~\ref{Fast_decay}).

The asymptotic behavior of discrete orthogonal polynomials was investigated in the last decades in numerous publications (cf. for example \cite{Shar1,Shar2,Shar3,Shar4,Shar5,Shar6,Nur}; a lot of important relations to many other domains of mathematics outside of signal processing and computational mathematics was found (cf. \cite{Apt-As, BKMM-2003} and the recent monograph \cite{BKMM-2007}). On the other hand, the results in these publications are limited to the asymptotic behavior of discrete orthogonal polynomials of {relatively small degree $n$} compared to the number $N$ of the grid points: $n < \alpha \sqrt{N}$. 
The fact of very rapid decay of their values near the endpoints for $n$ close enough to $N$ was neither studied nor even mentioned before. At the same time, this fact is extremely important for the analysis of the least squares polynomial approximations on discrete grids (with applications in signal processing), as it is shown in Section~\ref{sec-lsapprox-IGS}.

Numerical computations for the construction of the discrete orthogonal polynomials, used in the procedure of outlier detection, should be done with special care to numerical stability. We address this issue in Appendix~\ref{App-Hahn-computation}.


But first we explain why one shall use discrete orthogonal polynomials for removal of ``polynomial trends'' of sufficiently high degrees.

\section{Least squares polynomial fitting of  high degree: why discrete orthogonal polynomials?}


Removing linear and quadratic trends for more detailed study of smaller variations of time series of experimental data is very common in data analysis. 
This is usually done by finding the coefficients of the respective linear or quadratic polynomial that best fits the data using the least squares approximation techniques. 
Removing higher-order trends is more difficult. 
This section explains the theoretical and practical implications of higher-order polynomial fitting for the newcomers in the field of numerical analysis. 

First we should understand 
why is it not enough to use the least squares approximation directly finding the coefficients of some high-degree polynomial $ p (x) = a_0 + a_1 x + \ldots + a_n x^n $ of best approximation? 
The answer is very simple: the least squares approximation becomes \emph{very} numerically unstable resulting in a wrong set of the coefficients $a_i$ if one uses computations with double precision or even quadruple precision. Speaking in terms of numerical linear algebra this is reflected in the extremely large condition number of the corresponding least squares approximation matrix.

The reason for this can be easily explained geometrically. Take for example monomials of degrees 1, 10 and 30 on the interval $[0,1]$: \\
 \includegraphics[width = 0.3 \textwidth]{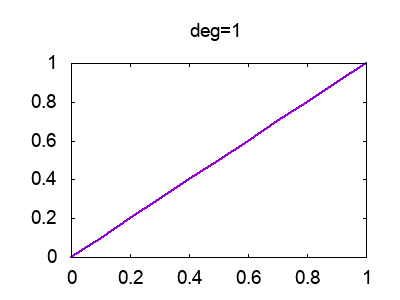}
 \includegraphics[width = 0.3 \textwidth]{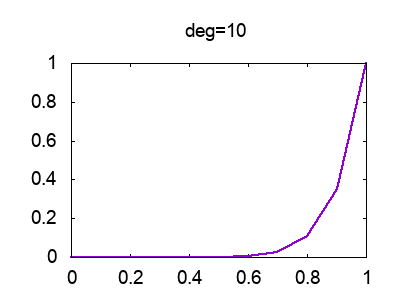}
 \includegraphics[width = 0.3 \textwidth]{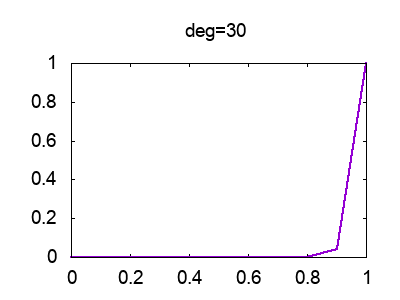}

Taken as vectors of their values at the points of the 
lattice $\{0, 0.1, 0.2, \ldots , 0.9, 1\}$
we easily understand that the second and the third vector are very close in the 11-dimensional space: their scalar product and their lengths are close to 1 so the angle between these two vectors is close to zero.
The same is obviously true for vectors representing monomials of higher degrees.

So the basis of vectors formed by such monomials geometrically is very bad for decomposition of an arbitrary data vector w.r.t.\ their basis. 
The matrix formed with such basis vectors (the Vandermond matrix in this monomial case) is very bad for inversion and other numerical linear algebra operations: 
the condition number of the Vandermond matrix for monomials of orders $0 \ldots 10$ is
$ C (W {}_{11}^{10})> 10^{8} $ and for orders $0 \ldots 30$ even
$ C (W {}_{31}^{30})> 10^{19} $!

It is obviously much more favorable to use orthogonal bases for the decomposition of data vectors. This is where discrete orthogonal polynomials appear: formally taking the basis of monomials (treated as finite dimensional vectors of their values on some finite grid of points) and using some orthogonalization algorithm, for example the Gram-Schmidt algorithm, one obtains a set of orthogonal vectors that can be viewed as the values of some polynomials on the same discrete grid --- exactly the set of so called Hahn polynomials (with parameters \(\alpha=\beta=0\), cf.\ \cite{NSU85,OLB10}). Unfortunately, instability of the least squares approximation algorithms for high degree polynomials is reflected in instability of the standard algorithms for the construction of discrete orthogonal polynomials of high degrees.

We address the problem of robust construction of Hahn polynomials in Appendix~\ref{App-Hahn-computation}.

In Section~\ref{Fast_decay}, we address some new theoretical properties of Hahn polynomials important for the interpretation of the results of best least squares polynomials approximation.

Before going into theoretical details, we expose in the next section the experimental results illustrating the effects that required these theoretical studies --- separation of ``slow'' dynamics of a signal (in our case, any of the coordinates of a GPS satellite with amplitude of tens of thousands of kilometers) from small (a few centimeters only) but ``fast'' variations of two types: outliers (isolated deviations) or ``jumps'' (step-like changes)  using the complete system of discrete orthogonal Hahn polynomials on equidistant grids as a filter.


\section{Least squares approximation of IGS orbits with discrete orthogonal polynomials}
\label{sec-lsapprox-IGS}


Analysing ``slow'' (polynomial) trends of a given time series
we use a reshaping of the approximation polynomial $p_m(t)$ of degree $m$ as a sum of normalized Hahn polynomials $\Q_{k}^N$ (cf.~section~\ref{DOPgeneral}) of degrees $k=0,\ldots, m$ orthogonal on the $(N+1)$-point grid where the analyzed time series $f_j = f(t_j)$ is given:
\begin{equation}\label{eq-pm}
    f(t_j) \approx p_m(t_j) = \sum_{k=0}^m b_k \Q_{k}^N(t_j).
\end{equation}
In such a case the coefficients $b_k$ of the polynomial of the best least squares approximation is readily found as scalar products on the data grid $t_j$, $j=0,\ldots, n$:
\begin{equation}\label{eq-bk}
 b_k = \sum_{j=0}^N f(t_j)\Q_{k}^N(t_j).
\end{equation}
This is a simple and robust computation, provided the correct values of $\Q_{k}^N$ at the grid points $t_j$ are pre-computed for the points of the (equidistant) grid.

In the computational experiments we analyzed one of the coordinates $X(t_j)$ of a GPS satellite (taken from an SP3 files published by the IGS \cite{1}):

\begin{minipage}{0.35 \textwidth}
    \includegraphics[width = 1.0 \textwidth]{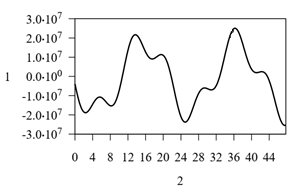}
  \end{minipage}
\begin{minipage}{0.6 \textwidth}
  {\small
    --- A typical graph of one of the coordinates of a GPS satellite (terrestrial rotating Cartesian coordinate system).
    \\
    Horizontal axis: time (in hours).
    \\
    Vertical axis: coordinate (in meters).
    }
\end{minipage}

Let us approximate this rapidly (but smoothly) varying time series by a polynomial of some high degree in a moving window of 4 consecutive days. Since the time step of the data in SP3 files is equal to 15 minutes, we have 96 data points for each day, or 384 points for 4 days. Constructing the orthogonal Hahn polynomials for this regular grid and finding the scalar products \eqref{eq-bk} we find the least squares polynomial approximations \eqref{eq-pm} of degrees up to 383. Then, taking the residues
\begin{equation}\label{eq-res}
res_m(t_j) = X(t_j) - p_m(t_j)
\end{equation}
we find the ``detrended'' results. As our experiments have shown, interesting results are obtained when we find the residues $res_m(t_j)$
for approximations by polynomials of degrees at least 175.

\begin{figure}[ht]
  \centering
  \includegraphics[width=0.9\textwidth]{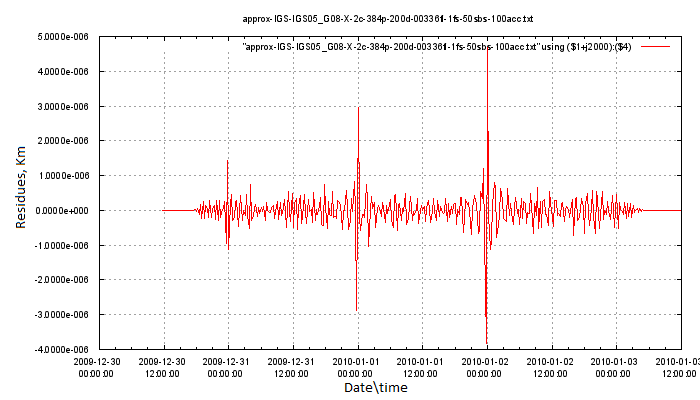}
  \caption{Residue for 4 days approximation of G08 orbit (30.12.2009--03.01.2010).}
  \label{Fig-1}
  \end{figure}

On Figure~\ref{Fig-1} given below (as a typical example) we plot the difference (the residue $res_{200}(t_j)$) between the approximating polynomial of the best least squares approximation of degree 200 and the X-coordinate of the orbit of the G08 satellite (GPS navigation system) on the interval from 12:00:00 on 30.12.2009 until 12:00:00 03.01.2010 (4 days) at 384 points from the SP3-files for these four days.

The vertical axis gives the residue (deviation of the orbit from  the approximating polynomial) in kilometers. So we see
\textbf{anomalous ``jumps'' at 00:00:00  with magnitude
$\sim $ 1 cm every day}, with chaotic oscillations with typical magnitudes 0.1--0.5~mm between them. Also we observe the important phenomenon of vanishing of the residue near the ends of the 4-day interval. 
The last phenomenon is due to the newly found effect of rapid decay of discrete orthogonal polynomials of high degree near the both ends of the grid, proved in Section~\ref{Fast_decay}. 
Indeed, since the systems of Hahn polynomials $\Q_{k}^N(t_j)$, $k=0,\ldots, N$, gives a complete orthogonal basis in the finite-dimensional space of all (scalar) data on the $(N+1)$-point grid,
the residue \eqref{eq-res} is the sum of the highest-degree Hahn polynomials:
\begin{equation}\label{eq-res1}
res_m(t_j) = \sum_{k=m+1}^{N} b_k \Q_{k}^N(t_j),
\end{equation}
rapid decay of the values of $\Q_{k}^N(t_j)$ of degree close to $N$ guarantees therefore the observed very rapid decay of the residue $res_{200}(t_j)$ near the boundaries of the 4-day time window.

The chaotic 0.5-mm oscillation of the graph of $res_{200}(t_j)$ between the day boundaries is explained by the method of formation of the SP3 orbit data: the trajectories are calculated numerically with the usual double-precision (giving approx. 15 decimal digits) and then the result is rounded to 1-mm precision (11 decimal digits as a maximum).

The nature of large ``jumps'' at the day boundaries is also easily explained. 
It is related to the processing methodology used by the IGS Analysis Centers: they collect measurements from the large ground network of observation stations for each day separately, then using rather sophisticated models
calculate various IGS data products \cite{1} for each day separately. The RMS precision of the orbit final data (as stated in \cite{1})
is  $\sim3$~cm, so ``jump patterns'' for the residues  $res_{200}(t_j)$ at the day boundaries are the results of slightly different (a few centimeters of difference) computational approximations of the \emph{physical} satellite orbits for two consecutive days after the ground data of the observation stations are processed.  Obviously such small jumps are well within the stated RMS precision $\sim3$~cm.

In our case these ``jump patterns'' visible in the approximation residues for a smooth time series with a discontinuity can be easily explained by the following model approximation experiments. In fact this effect is one of the manifestations of the well-known Gibbs phenomenon (cf. for example \cite{Gibbs-phenomenon}) in approximation theory. 


\begin{center}
    \textbf{Discontinuities (jumps) in a data series}
    \end{center}
    
Let us take a ``jump'' of magnitude 1 at the point $t=40$ in a series of 101 points $t=0,\ldots, 100$ (the left Figure below):

\begin{minipage}{0.4\textwidth}
\includegraphics[width=0.9\textwidth]{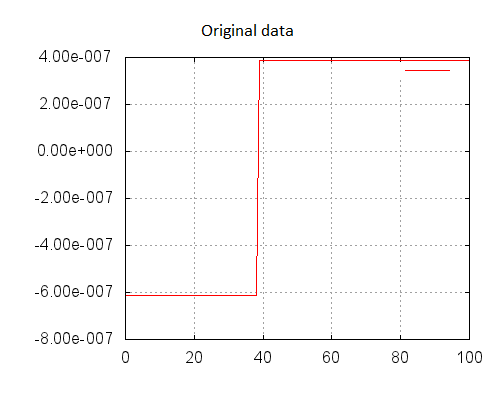}
\end{minipage}
\begin{minipage}{0.4\textwidth}
\includegraphics[width=1.2\textwidth]{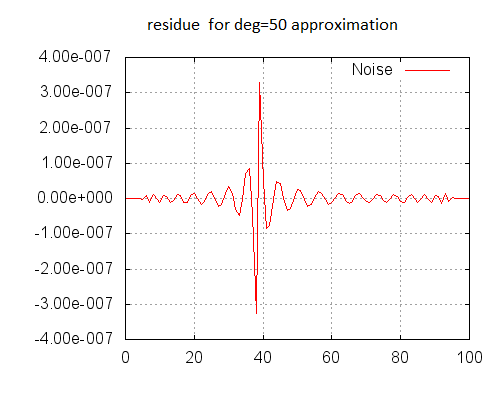}
\end{minipage}

After least squares approximation with a polynomial of degree 50, the residue $res_{50}(t_j)$,
as shown on the right picture above, has decaying oscillations around the jump point $t=40$,
and a specific pattern at this point --- two (almost-)symmetric spikes of magnitude 0.33
(a little less that the half of the magnitude 1 of the original jump). This "jump pattern''
is easily recognizable at the day boundaries on the previous Figure showing the residue $res_{200}(t_j)$ of
approximation of the G08 coordinate for the period of 4 days.

But such small variations at the day boundaries turn out to be not the only anomalies in the IGS orbit SP3 data.

On the next Figure~\ref{Fig-Huge-jump} given below we plot the residue $res_{200}(t_j)$ between the approximating polynomial of the best least squares approximation of degree 200 and the X-coordinate of the orbit of the G25 satellite on the interval from 12:00:00 on 31.07.2010 until 00:00:00 on 05.08.2010:

\begin{figure}[ht]
  \centering
  \includegraphics[width=0.9\textwidth]{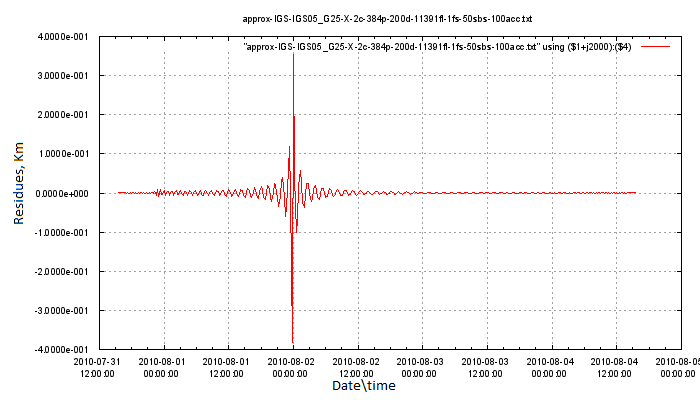}
  \caption{Huge anomalous ``jump'' for G25 at 00:00:00, 02.08.2010.}
  \label{Fig-Huge-jump}
  \end{figure}

Here we observe a huge anomalous ``jump'' at 00:00:00, 02.08.2010  with approximate magnitude
 0.8 km with the typical ``jump pattern''.

Using the same methodology one can detect also another residue pattern, corresponding to an outlier in a time series.
For this  we have to make another computational experiment.


\begin{center}
\textbf{Anomalous single value (outlier)  in a data series}
\end{center}

Let us take an ``outlier'' of magnitude 1 at   $t=40$ in a series of 101 points $t=0,\ldots, 100$ (the left Figure below):

\begin{minipage}{0.45\textwidth}
\includegraphics[width=0.9\textwidth]{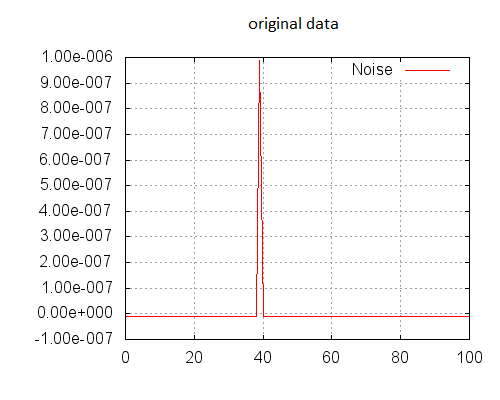}
\end{minipage}
\
\begin{minipage}{0.4\textwidth}
\includegraphics[width=1.2\textwidth]{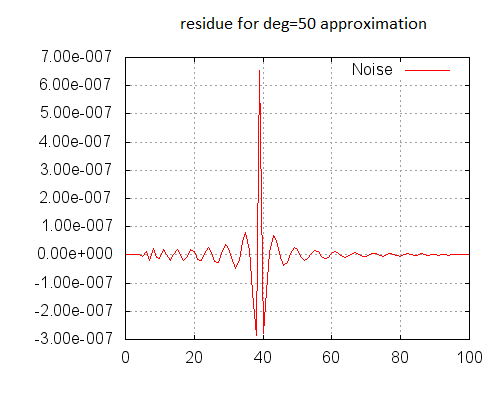}
\end{minipage}

After least squares approximation with a polynomial of degree 50, the residue $res_{50}(t_j)$,
as shown on the right picture above, has a pattern, different from the previous ``jump pattern''  around the point $t=40$,
--- one positive spike at $t=40$ and (almost-)symmetric negative smaller spikes; their total magnitude being
a bit less that the magnitude 1 of the original outlier.

This ``outlier pattern''
is also found (with rather large magnitudes \emph{not} at the day boundaries in the residues $res_{200}(t_j)$ after polynomial approximation of the GPS orbits published by the IGS for the year 2010.

They are very similar to either  ``outlier pattern'' or  ``jump pattern'' described above.

\textbf{List of various anomalies for 2010 of GPS final orbits
(\emph{not} at the day boundaries)}

The processed file names listed below describe anomalies for the $X$ coordinate with the ``outlier pattern'' (starts with ``impuls-'' in the file name), its magnitude in centimeters or meters (50cm for the first file), and the satellite number (G09 for the first file). Here we list only the large (more that 20 cm) anomalies found after processing of the SP3 data for the year 2010:

\noindent
 impuls-50cm-G09-X-2c-384p-200d-10721fl-1fs-50sbs-100acc  \\
 impuls-50cm-G27-X-2c-384p-200d-30821fl-1fs-50sbs-100acc  \\
 impuls-5m-G10-X-2c-384p-200d-22446fl-1fs-50sbs-100acc  \\
 impuls-5m-G20-X-2c-384p-200d-25796fl-1fs-50sbs-100acc  \\
 impuls-5m-G25-X-2c-384p-200d-10721fl-1fs-50sbs-100acc  \\
 impuls-5m-G30-X-2c-384p-200d-02011fl-1fs-50sbs-100acc  \\
 impuls-10m-G04-X-2c-384p-200d-33836fl-1fs-50sbs-100acc  \\
 impuls-10m-G22-X-2c-384p-200d-13066fl-1fs-50sbs-100acc  \\
 impuls-10m-G32-X-2c-384p-200d-07371fl-1fs-50sbs-100acc  \\
 impuls-15m-G01-X-2c-384p-200d-28476fl-1fs-50sbs-100acc  \\
 impuls-15m-G03-X-2c-384p-200d-23116fl-1fs-50sbs-100acc  \\
 impuls-15m-G16-X-2c-384p-200d-20436fl-1fs-50sbs-100acc  \\
 impuls-15m-G32-X-2c-384p-200d-33501fl-1fs-50sbs-100acc  \\
 impuls-30m-G02-X-2c-384p-200d-09716fl-1fs-50sbs-100acc  \\
 impuls-30m-G11-X-2c-384p-200d-27806fl-1fs-50sbs-100acc  \\
 impuls-30m-G17-X-2c-384p-200d-03351fl-1fs-50sbs-100acc  \\
 impuls-40m-G14-X-2c-384p-200d-22781fl-1fs-50sbs-100acc  \\
 impuls-50m-G29-X-2c-384p-200d-01676fl-1fs-50sbs-100acc  \\
 impuls-50m-G31-X-2c-384p-200d-31156fl-1fs-50sbs-100acc  \\
 impuls-100m-G01-X-2c-384p-200d-18426fl-1fs-50sbs-100acc  \\



\begin{figure}[ht]
  \centering
  \includegraphics[width=0.8\textwidth]{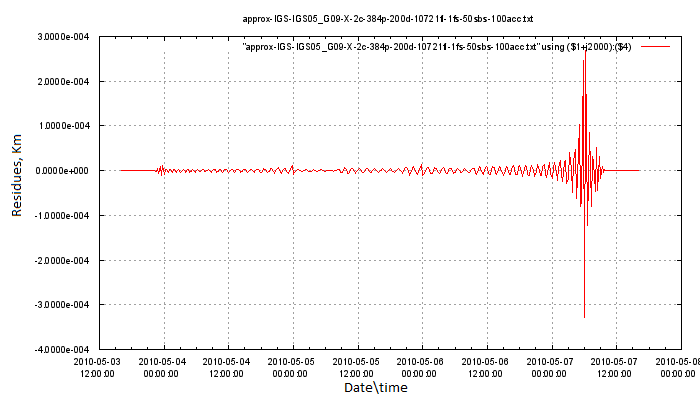}
  \caption{Satellite maneuver pattern.}
  \label{Fig-impulsN}
  \end{figure}

Also we give below on Figure~\ref{Fig-impulsN} a bit different anomaly with magnitude 50~cm for GPS satellite G09, time interval: 04.05.2010--07.05.2010, again visible on the graph of the residue $res_{200}(t_j)$.

Most (but not all!) of such large anomalies can be explained by maneuvering of the GPS satellites (the announcements of the planned maneuvers are published by the GPS Monitor and Control Segment). 




Summing up, we observed a few interesting phenomena apparent in approximation of presumably smooth time series of GPS orbits with high-degree polynomials:
\begin{itemize}
\item small ``digital noise'' of order 0.5~mm as the result of rounding of the computed orbits to the nearest millimeter value;
\item regular jumps of order $\sim1$~cm at the day boundaries (Gibbs phenomenon) due to the processing methodology of IGS Analysis Centers as well as other jumps and outliers not explainable by their processing routines; 
\item anomalously fast decay of the approximation residues near the end points of the approximation interval. This is a new numerical (not physical!) phenomenon explained in the next sections.
\end{itemize}


\section{Mathematical background and tools}
\label{sec-mathbt}

The approximation methodology used in the previous section was based on \emph{least squares approximation of IGS orbits with discrete orthogonal polynomials}. At the moment we choose equal weights at the discrete time grid.

There is extensive literature on polynomials, orthogonal on intervals or systems of intervals of the real line. Below we briefly recall the properties of such orthogonal polynomials, which are usually called continuous orthogonal polynomials and expose very important differences with  polynomials, orthogonal on finite subsets of points of the real line -- discrete orthogonal polynomials.

\subsection{Continuous orthogonal polynomials -- a brief reminder}
\label{ContOrthPol}


Legendre polynomials (orthogonal polynomials on the interval $[-1,1]$ with unit weight) are very well known (we mention here only two basic texts \cite{Szego, Askey}), their properties include:
\begin{itemize}
 \item they satisfy a 2nd-order recurrence;
 \item they have all their roots in the interval  $[-1,1]$;
 \item their local extrema on $[-1+\varepsilon,1-\varepsilon]$ (where $\varepsilon$ depends on the degree of a polynomial) are not very different in magnitude with higher absolute values closer to the endpoints.
\end{itemize}
For example, we plot the Legendre polynomial of degree 30 in the interval \([-1,1]\) (Figure~\ref{Fig-Legendre30}):


\begin{figure}[ht]
\centering
\includegraphics[width=0.9\textwidth]{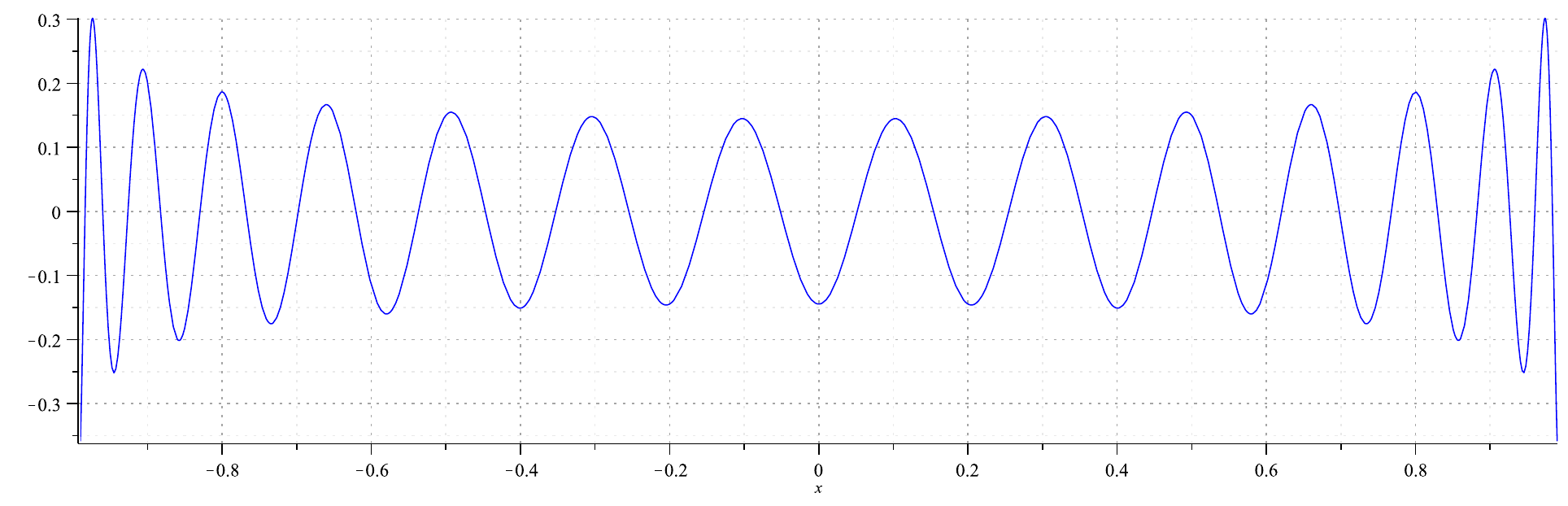}
\caption{Legendre polynomial of degree 30.}
\label{Fig-Legendre30}
\end{figure}

Continuous orthogonal polynomials are extensively used in theoretical and applied research and are ubiquitous in mathematical physics. On the contrary, in modern digital signal processing, image processing and many other computational applications data are usually given on finite discrete grids (either equally spaced or irregular grids). As a consequence, standard theoretical results (least squares approximations etc.) and algorithms should be adapted to finite-dimensional functional spaces on such grids, in particular spaces of polynomials with arguments on such finite grids --- discrete orthogonal polynomials. 


\subsection{Discrete orthogonal polynomials}
\label{DOPgeneral}


In the following, we study only  equidistant finite grids \(X=\{0,1,\ldots,N\}\) suitable for our case study of GPS and GLONASS orbits.

Discrete orthogonal polynomials are the polynomials $p_m(x) = \sum_{k=0}^{m} a_k x^k$ of degrees $m=0,\ldots,N$, with values taken on equidistant finite grid \(x_j\in X=\{0,1,\ldots,N\}\), orthogonal with respect to the standard inner product (with weights $w_j$):
\begin{equation}\label{eq-wei}
  \sum_{j=1}^n p_m(x_j)p_s(x_j) w_j = 0, \quad m \neq s .
\end{equation}

Formally, many of their properties are similar to those of continuous orthogonal polynomials:
\begin{itemize}
 \item they satisfy a 2nd-order recurrence;
 \item they have all their roots in the interval  $[0,N]$.
\end{itemize}
Details can be found in \cite{NSU85, BKMM-2007}.

In this paper, we study some properties of Hahn polynomials (also known as Chebyshev discrete polynomials cf.~\cite{NSU85}) with equal weights $w_j\equiv 1$.




For convenience, we generalize the notation of the Pochhammer symbol. Let
\begin{equation}\label{a_k}
(a)_{k}=\begin{cases} a(a+1) \ldots (a+k-1), & \text{when}\ k>0, \\ 1, & \text{when}\ k=0, \\ a(a-1) \ldots (a+k+1), & \text{when}\ k<0. \end{cases}
\end{equation}
Clearly, 
\begin{equation}\label{-a_k}
(-a)_{k}=(-1)^{k}(a)_{-k}.
\end{equation}

Let \(N\in\mathbb{Z}_{+}\) be fixed, and \(X=\{0,1,\ldots,N\}\). For \(n,x\in X\), the Hahn (discrete Chebyshev) polynomial \(Q_{n}^{N}(x)\) of degree \(n\) \emph{with equal weights $w_j\equiv 1$} is defined as
\begin{equation}\label{Qn}
Q_{n}^{N}(x)=\sum_{k=0}^{n}\frac{(-1)^{k}(n)_{-k}(n+1)_{k}(x)_{-k}}{(k!)^2(N)_{-k}}.
\end{equation}

\begin{remark}\label{Qn-rem}
Although orthogonal polynomials can be defined in terms of hypergeometric functions, the authors believe, that this way may contain some redundancies that are not essential for conveying the ideas of this work. Nevertheless, the reader may wish to check (which should be pretty straightforward) that \eqref{Qn} coincides with Formula~(18.20.5) from \cite{OLB10} for \(\alpha=\beta=0\).
\end{remark}

The family \(\left\{Q_{n}^{N}(x)\right\}_{n\in X}\) is orthogonal with the weight function \(w_{x}\equiv 1\) and with the square of the norm given by
\begin{equation}\label{Qn-norm}
h_{n}^{N}=\left\|Q_{n}^{N}(x)\right\|^2=\sum_{x\in X}\left(Q_{n}^{N}(x)\right)^2=\frac{(N+1)_{n+1}}{(2n+1)(N)_{-n}}.
\end{equation}

Below we use the following notation for the \emph{normalized Hahn polynomials}:
\begin{equation}\label{Qn-normalized}
  \Q_{n}^{N}(x)= Q_{n}^{N}(x) \big/ \sqrt{h_{n}^{N}}.
  \end{equation}

For example, in Figure~\ref{Fig-Hahn30} we plot (normalized) Hahn polynomial \(\Q_{30}^{30}(x)\) of degree~30 on \(X=\{0,1,\ldots,30\}\) (31 equidistant points), the values of the polynomial between the points of the grid are omitted, the values of the polynomial at the grid points are connected by straight line segments:



\begin{figure}[ht]
\centering
\includegraphics[width=0.9\textwidth]{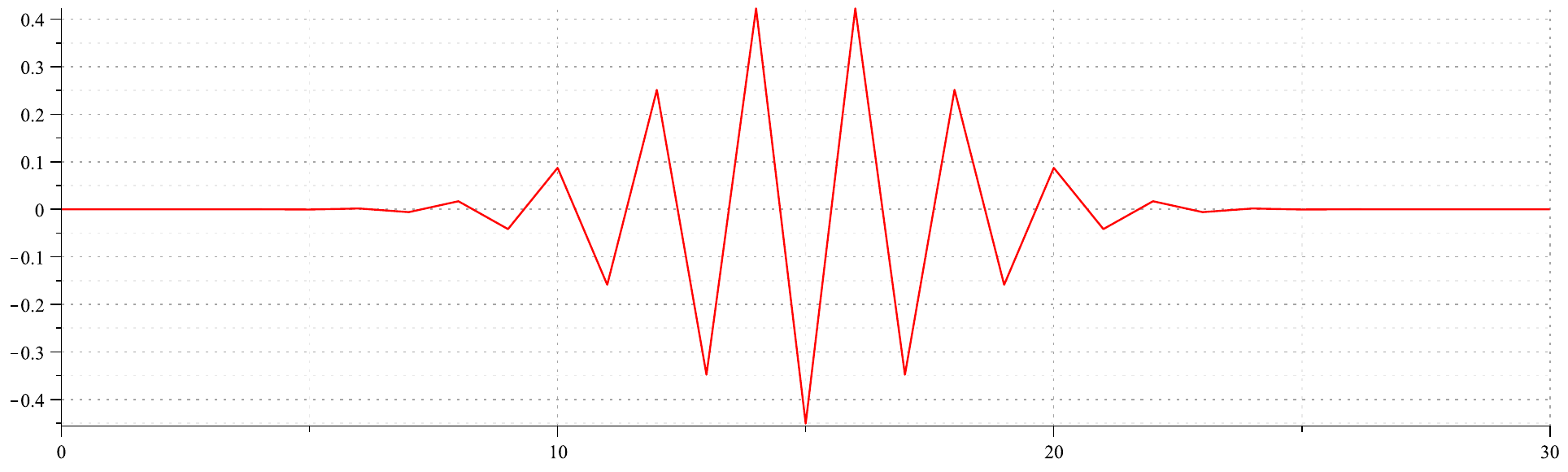}
\caption{Normalized Hahn polynomial \(\Q_{30}^{30}(x)\) of degree 30.}
\label{Fig-Hahn30}
\end{figure}

As one can see, contrary to the continuous Legendre polynomials, the values of the normalized Hahn polynomials \emph{at the points of the grid} almost vanish near the boundary, here are some values near $x=0$:
\vspace{-5mm}
\begin{center}
\begin{tabular}[t]{|c|c|c|c|c|}
\hline
\(x\) & 0 & 1 & 2 & 3 \\ \hline \(\Q_{30}^{30}(x)\) & \(2.9079 \cdot 10^{-9}\) & \(-8.7236 \cdot 10^{-8}\) & \(1.2649 \cdot 10^{-6}\) & \(-1.1806 \cdot 10^{-5}\) \\ \hline
\end{tabular}
\end{center}

\begin{figure}[ht]
\centering
\includegraphics[width=0.9\textwidth]{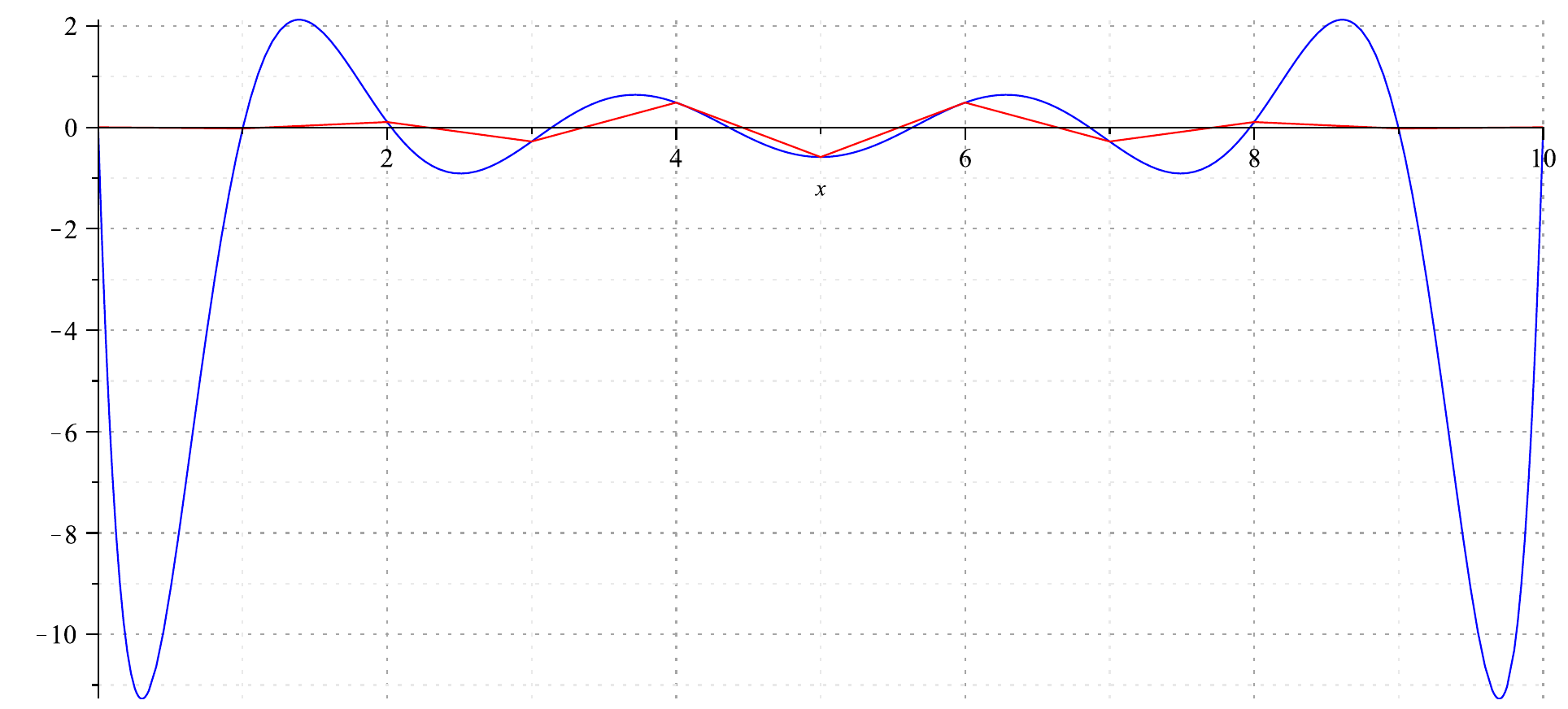}
\caption{Values of normalized Hahn polynomial \(\Q_{10}^{10}(x)\) of degree 10 between integer points.}
\label{Fig-Hahn10cont}
\end{figure}

The values of this polynomial between the grid points are, on the contrary, very large. In the table below, we give some values of \(\Q_{30}^{30}(x)\) at \(x=\frac{1}{2}+k\) for \(k=0,1,2,3\):
\vspace{-5mm}
\begin{center}
\begin{tabular}[t]{|c|c|c|c|c|}
\hline
\(x\) & 0.5 & 1.5 & 2.5 & 3.5 \\ \hline \(\Q_{30}^{30}(x)\) & \(-1.2398 \cdot 10^{6}\) & \(67920.4\) & \(-6460.054\) & \(898.31\) \\ \hline
\end{tabular}
\end{center}
The fact is also illustrated on Figure~\ref{Fig-Hahn10cont} where we plot (as the blue curve) the polynomial \(\Q_{10}^{10}(x)\) on the continuous interval $x=[0,10]$, again (as in the previous plot) connecting the values of the polynomial at the grid points by red straight line segments. As we can see, close to the boundary the polynomial \(\Q_{10}^{10}(x)\) oscillates between the grid points with very large amplitude while its values in the integer points are very small.

One can observe the same pattern  for \(\Q_{n}^{N}(x)\) with sufficienly large $N$ and the degree $n$ close to $N$: the values at integer points vanish close to the endpoints of the interval while the values in between take rather high values. We finalize this observation with the table of the orders of values of \(\Q_{75}^{100}(x)\) (the order \(10^{k}\) means that the precise value is \(\gamma \cdot 10^{k}\) with \(1\leqslant |\gamma| < 10\).)
\vspace{-5mm}
\begin{center}
\begin{tabular}[t]{|c|c|c|c|c|c|c|c|c|c|c|c|}
\hline
\(x\) & 0 & 0.5 & 1 & 1.5 & 1 & 2.5 & 3 & 3.5 & 4 & 4.5 & 5 \\ \hline \(\Q_{75}^{100}(x)\) & \(10^{-14}\) & \(10^{10}\) & \(10^{-12}\) & \(10^{9}\) & \(10^{-11}\) & \(10^{7}\) & \(10^{-10}\) & \(10^{6}\) & \(10^{-9}\) & \(10^{5}\) & \(10^{-8}\) \\ \hline
\end{tabular}
\end{center}

\begin{figure}[ht]
\centering
\includegraphics[width=0.9\textwidth]{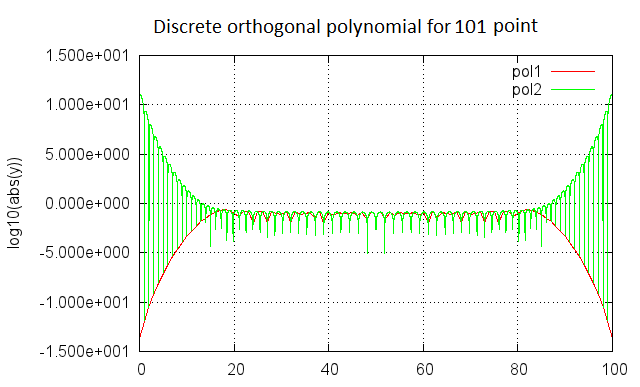}
\caption{Values of normalized Hahn polynomial \(\Q_{75}^{100}(x)\) of degree 75 for 101 points on the \(\log_{10}\) scale.}
\label{Fig-Hahn100-75-log}
\end{figure}

A much better understanding of this kind of behavior is obtained when one plots the same discrete orthogonal polynomials in \emph{log-scale}. For example, for the equidistant grid with 101 point and the normalized Hahn  (Chebyshev) polynomial polynomial of \(\deg=75\) the log-plot of \(\log_{10} |\Q_{75}^{100}(x)|\) is given in Figure~\ref{Fig-Hahn100-75-log}. Here the graph of the logarithm of the absolute value of the polynomial at integer points is given in red and between the integer points is in green color.

\begin{figure}[ht]
\centering
\includegraphics[width=0.9\textwidth]{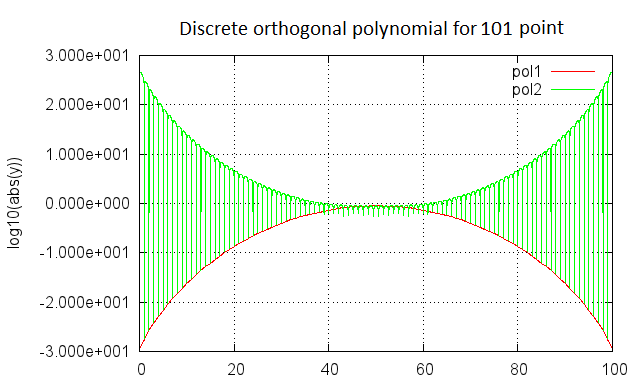}
\caption{Values of normalized Hahn polynomial of degree 100 for 101 points on the \(\log_{10}\) scale.}
\label{Fig-Hahn100-100-log-png}
\end{figure}

For the same grid with 101 points the log-plot of the Hahn (Chebyshev) polynomial of the maximal possible \(\deg=100\) is given in Figure~\ref{Fig-Hahn100-100-log-png}. Its values at the grid points near the boundary are as small as $10^{-30}$ and it oscillates between the integer points with the magnitude as big as $10^{30}$, approximately.

We now proceed to the formal proofs of these facts.


\section{Rapid decay of discrete orthogonal polynomials}
\label{Fast_decay}


We first give an idea which explains the mentioned above behavior of Hahn polynomials.

Note that the value \(Q_{n}^{N}(x)\) at \(x=m\) contains only \(\min\{m,n\}+1\) first summands, namely
\begin{equation}\label{Qn(m)}
Q_{n}^{N}(m)=\sum_{k=0}^{\min\{m,n\}}\frac{(-1)^{k}(n)_{-k}(n+1)_{k}(m)_{-k}}{(k!)^2(N)_{-k}},
\end{equation}
since the products \((n)_{-k}(m)_{-k}\) vanish for \(k\geqslant\min\{m,n\}+1\) (see the notation \eqref{a_k}).
Since for our purposes (in signal processing applications) \(n\) is relatively large, below we consider, for simplicity, \(m\leqslant n\).

Write down \(Q_{n}^{N}(m)\) for a few smallest values of \(m\):
\[
Q_{n}^{N}(0)=1, \quad Q_{n}^{N}(1)=1-\frac{n(n+1)}{N}, \quad Q_{n}^{N}(2)=1-\frac{2n(n+1)}{N}+\frac{(n-1)n(n+1)(n+2)}{2N(N-1)}.
\]
We see that \(Q_{n}^{N}(0)\) has only the first summand corresponding to \(k=0\). The rest vanish since \((0)_{-k}=0\) for \(k\geqslant 1\). Similarly, \(Q_{n}^{N}(1)\) has only the summands corresponding to \(k=0\) and \(k=1\) because \((1)_{-k}=0\) for \(k\geqslant 2\). In general, the following estimates hold:
\[
    \left|\Q_{n}^{N}(0)\right|=
\left|\frac{Q_{n}^{N}(0)}{\sqrt{h_{n}^{N}}}\right|=\frac{1}{\sqrt{h_{n}^{N}}}, \quad 
\left|\Q_{n}^{N}(1)\right|=
\left|\frac{Q_{n}^{N}(1)}{\sqrt{h_{n}^{N}}}\right|\leqslant \frac{N+2}{\sqrt{h_{n}^{N}}},\ldots
\]
From this observation, we can get an idea why the values of the normalized Hahn polynomials are small at the integer points close to the endpoints.

At the same time, if we consider the values of \(Q_{n}^{N}(x)\) in between the integer points, then the situation is quite different. Now the terms \((m)_{-k}\) and \((n)_{-k}\) does not vanish. Moreover, it is easy to see, that, for a fixed \(k\), the term \((m)_{-k}\) takes the largest value when \(m\in(0,1)\), for instance \((0.3)_{-10}\approx -42884\). At the same time, one can estimate, that the value
\[
\frac{(n)_{-k}(n+1)_{k}}{(k!)^2(N)_{-k}}>1
\]
for \(n\geqslant \frac{N}{2}\) (in the applications, we consider \(n\approx \frac{3}{4}N\)) and small values of \(m\): \(k\leqslant m\leqslant \frac{N}{10}\).

We now proceed with estimating the values \(\Q_{n}^{N}(m)\) near the endpoints.

For fixed \(N\), \(n\) and \(m\), consider the absolute value of the \(k\)-th summand of \eqref{Qn(m)} as a function of~\(k\):
\[
q_{n, m}^{N}(k)=\frac{(n)_{-k}(n+1)_{k}(m)_{-k}}{(k!)^2(N)_{-k}}, \quad k=0,\ldots,\min\{m,n\}.
\]

\begin{figure}[ht]
\centering
\includegraphics[width=0.9\textwidth]{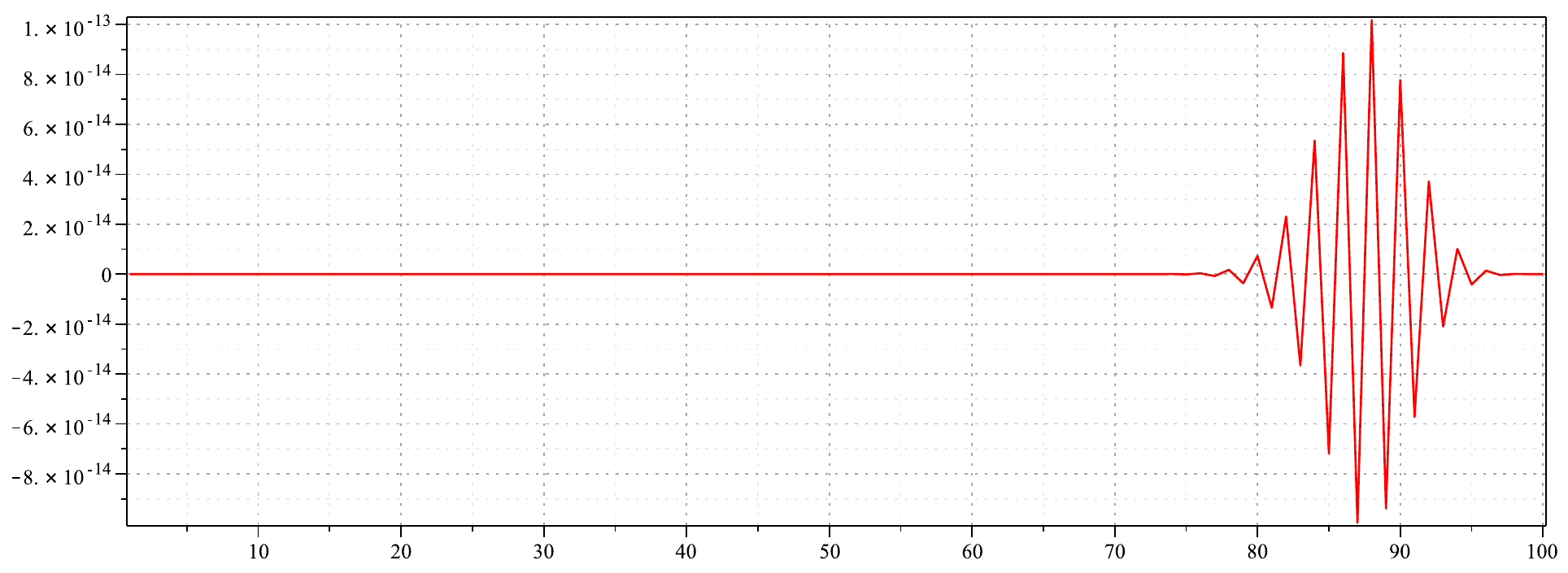}
\caption{Values of \(q_{750, 100}^{1000}(k)\).}
\label{Fig-Q_k_1000-750-100}
\end{figure}

Let the maximum value of \(q_{n, m}^{N}(k)\) be attained at the point \(\widetilde{k}\) and be equal to 
\[
\widetilde{q}_{n, m}^{N}=|q_{n, m}^{N}(\widetilde{k})|=\max_{k}{|q_{n, m}^{N}(k)|}.
\]
There is a simple algorithm for finding \(\widetilde{k}\). Indeed, for \(k=1,\ldots,\min\{m,n\}\), denote
\begin{equation}\label{HanhKratio}
r_{n, m}^{N}(k)=\left|\frac{q_{n, m}^{N}(k)}{q_{n, m}^{N}(k-1)}\right|=\frac{(n-k+1)(n+k)(m-k+1)}{k^2(N-k+1)}.
\end{equation}
We look for such values of \(k\) that \(r_{n, m}^{N}(k)\leqslant 1\) which means that \(\left|q_{n, m}^{N}(k)\right|\leqslant \left|q_{n, m}^{N}(k-1)\right|\). Since the denominator of \eqref{HanhKratio} is always positive,  the equation \(r_{n, m}^{N}(k)= 1\) is equivalent to
\[
(n-k+1)(n+k)(m-k+1)-k^2(N-k+1)=0
\]
or \(f(k)=0\) where
\begin{equation}\label{f(k)}
f(k)=2k^3-(N+m+3)k^2+\bigl((m+1)-n(n+1)\bigr)k+(m+1)n(n+1).
\end{equation}

Note that \(f(0)=(m+1)n(n+1)\geqslant 0\) for positive \(n\). On the other hand, 
\[
f(m+1)=(m+1)^2(m-N)\leqslant 0
\]
which yields the existence of a root of \(f(k)\) on \([0,m+1]\).


We now study \(f(k)\) for the typical values of the parameters \(n \sim 0.75N\), \(m \sim 0.1N\) (i.e., large enough \(n\) and small enough \(m\)). Substitute \(\tilde{N}=N+2\), \(\tilde{n}=n(n+1)\), \(\tilde{m}=m+1\) in \eqref{f(k)}. Then
\[
f(k)=2k^3-(\tilde{N}+\tilde{m})k^2+(\tilde{m}-\tilde{n})k+\tilde{m}\tilde{n}.
\]
If \(n\geqslant 1\), then the fact that \(f(-\infty)<0\) and \(f(0)=\tilde{m}\tilde{n}>0\) yields the existence of a real root \(k_{(1)}<0\). Furthermore, 
\[
f(\tilde{m})=\tilde{m}^3-\tilde{N}\tilde{m}^2+\tilde{m}^2=(\tilde{m}+1-\tilde{N})\tilde{m}^2<0
\]
for \(\tilde{m} \sim 0.1\tilde{N}\). This shows that \(f(k)\) has another real root in the interval \((0,\tilde{m})\). We now note that for \(\tilde{m} \sim 0.1\tilde{N}\), the first and the second derivatives,
\[
f'(k)=6k^2-2(\tilde{N}+\tilde{m})k+(\tilde{m}-\tilde{n})
\]
and
\[
f''(k)=12k-2(\tilde{N}+\tilde{m})=2(6k-\tilde{N}-\tilde{m}),
\]
are negative on \((0,\tilde{m})\) so that \(f(k)\) is decreasing and concave down on this interval. Consequently, to estimate the root \(k_{(2)}\in (0,\tilde{m})\), we can use the secant method (with the initial point \(k_{0}=\tilde{m}\) and the first iteration \(k_{1}=0\)) for estimating \(k_{(2)}\) from below:
\[
k_{2}=k_{1}-\frac{f(k_{1})(k_{1}-k_{0})}{f(k_{1})-f(k_{0})}=0-\frac{f(0)(0-\tilde{m})}{f(0)-f(\tilde{m})}=\frac{\tilde{m}^2\tilde{n}}{\tilde{m}\tilde{n}-(\tilde{m}+1-\tilde{N})\tilde{m}^2}=\frac{\tilde{m}\tilde{n}}{\tilde{n}-(\tilde{m}+1-\tilde{N})\tilde{m}}.
\]
For estimating \(k_{(2)}\) from above, we use Newton's method with the initial point \(k_{0}=\tilde{m}\):
\[
k_{1}=k_{0}-\frac{f(k_{0})}{f'(k_{0})}=\tilde{m}-\frac{(\tilde{m}+1-\tilde{N})\tilde{m}^2}{6\tilde{m}^2-2(\tilde{N}+\tilde{m})\tilde{m}+\tilde{m}-\tilde{n}}=\frac{3\tilde{m}^3-\tilde{N}\tilde{m}^2-\tilde{m}^2+\tilde{m}-\tilde{n}}{4\tilde{m}^2-2\tilde{N}\tilde{m}+\tilde{m}-\tilde{n}}
\]

Then, clearly, the estimate
\[
\left|Q_{n}^{N}(m)\right| \leqslant m\cdot\widetilde{q}_{n, m}^{N}
\]
holds for \(m\geqslant 1\).

In the table below, we list the mentioned values for \(N=100\), \(n=75\).
\vspace{-5mm}
\begin{center}
\begin{tabular}[t]{|c|c|c|c|c|c|c|}
\hline
\(m\) & 1 & 2 & 3 & 4 & 5 \\ 
\hline \(\widetilde{k}\) & 1 & 2 & 3 & 4 & 5 \\ 
\hline \(\left|Q_{n}^{N}(m)\right|\) & \(1.16\cdot 10^{-12}\) & \(3.2\cdot 10^{-11}\) & \(5.6\cdot 10^{-10}\) & \(7\cdot 10^{-9}\) & \(7\cdot 10^{-8}\) \\ 
\hline \(m\cdot\widetilde{q}_{n, m}^{N}\) & \(1.18\cdot 10^{-12}\) & \(6.8\cdot 10^{-11}\) & \(2\cdot 10^{-9}\) & \(4\cdot 10^{-8}\) & \(6\cdot 10^{-7}\) \\ 
\hline \hline
\(m\) & 6 & 7 & 8 & 9 & 10 \\ \hline \(\widetilde{k}\) & 6 & 7 & 7 & 8 & 9 \\ 
\hline \(\left|Q_{n}^{N}(m)\right|\) & \(6\cdot 10^{-7}\) & \(4\cdot 10^{-6}\) & \(7\cdot 10^{-5}\) & \(1\cdot 10^{-4}\) & \(4\cdot 10^{-4}\) \\ 
\hline \(m\cdot\widetilde{q}_{n, m}^{N}\) & \(7\cdot 10^{-6}\) & \(7\cdot 10^{-5}\) & \(6\cdot 10^{-4}\) & \(6\cdot 10^{-3}\) & \(5\cdot 10^{-2}\) \\ 
\hline
\end{tabular}
\end{center}
The table shows, in particular, that the estimate \(m\cdot\widetilde{q}_{n, m}^{N}\) is quite accurate for the small values of \(m\).



\section{Conclusion}

The above numerical experiments demonstrated the possibility of recognition of very small outliers and jumps (with relative magnitude $10^{-9}$ in comparison to the typical data values). 
Some of the results of our numerical experiments with GPS and GLONASS ephemeris day boundary discontinuities were briefly exposed in~\cite{Lobanov}.

An important new phenomenon of rapid decay of the residues of the least squares approximations by polynomials of high degree near the boundaries was observed; this phenomenon masks possible outliers near the boundaries of the processing window. One can prove appropriate estimates for the values of Hahn polynomials between the lattice points as well as estimates on the positions of their roots (to be reported in a later publication).

A numerically stable method for construction of discrete orthogonal polynomials of high degrees is proposed (Appendix~\ref{App-Hahn-computation}).


\begin{thebibliography}{9}
\bibitem{1}  \url{http://www.igs.org/products}

\bibitem{GriffithsRay2009}   \emph{Griffiths J., Ray J.R.} On the precision and accuracy of IGS orbits. Journal of Geodesy. 2009. V. 83, P.~277--287. DOI:10.1007/s00190-008-0237-6.

\bibitem{Lobanov} \emph{S.P.Tsarev, S.A.Lobanov}, Discontinuities and anomalous values in IGS final orbits,
Proc. Reshetnev Conf., Krasnoyarsk, 2013. v.~2. p.~120--122.

\bibitem{Szego} \emph{Szego G.} Orthogonal polynomials, In: AMS Colloquium Publi, vol. 79. 1939.

\bibitem{Askey} \emph{Askey, R.} Orthogonal polynomials and special functions. SIAM, 1975.

\bibitem{NSU85} \emph{A.F. Nikiforov, V.B. Uvarov, S.K. Suslov},
\newblock  Classical orthogonal polynomials of a discrete variable, Springer, 1991, 215~p.

\bibitem{Shar1} \emph{Sharapudinov I.I.} Some properties of polynomials, orthogonal on a finite system of points,
Izv. Vyssh. Uchebn. Zaved. Mat., 1983, no. 5, p.~85--88.

\bibitem{Shar2} \emph{Sharapudinov I.I.} On the best approximation and polynomials of the least quadratic deviation.
Analysis Mathematica. 1983, No~9(3), p.~223--234.

\bibitem{Shar3} \emph{Sharapudinov I.I.} Asymptotic properties of orthogonal Hahn polynomials in a discrete variable,
Mat. Sbornik, 180:9 (1989), p.~1259--1277.


\bibitem{Shar4} \emph{Sharapudinov I.I.} Approximation of discrete functions, and Chebyshev polynomials orthogonal on a uniform grid,
 Mat. Zametki, 67:3 (2000), p.~460--470.

\bibitem{Shar5} \emph{Sharapudinov I.I.} On the asymptotic behavior of Chebyshev orthogonal polynomials of a discrete variable,
Mat. Zametki, 48:6 (1990), p.~150--152

\bibitem{Shar6} \emph{Sharapudinov, I.I.} Convergence of the method of least squares.
 Mathematical Notes v.~53, no.~3 (1993), p.~335--344.

\bibitem{Nur}  \emph{Nurmagomedov  A.A.} About asymptotics of polynomials, orthogonal on arbitrary grids,
Izv. Saratov Univ. (N.S.), Ser. Math. Mech. Inform., 8:1 (2008), p.~25--31.

\bibitem{Apt-As}  \emph{Aptekarev, A.I., Van Assche, W.} Asymptotics of discrete orthogonal polynomials and the continuum limit of the Toda lattice. Journal of Physics A: Mathematical and General, 2001, v.~34(48), No.10627.

\bibitem{BKMM-2007} \emph{Baik, J., Kriecherbauer, T., McLaughlin, K. D. R., Miller, P. D.}
Discrete Orthogonal Polynomials. Princeton University Press, 2007, 170~p.

\bibitem{BKMM-2003} \emph{Baik, J., Kriecherbauer, T., McLaughlin, K.R. Miller, P.D.} Uniform asymptotics for polynomials orthogonal with respect to a general class of discrete weights and universality results for associated ensembles: announcement of results. International Mathematics Research Notices, 2003(15), p.821--858.

\bibitem{OLB10} \emph{Olver,~F.~W.~J., Lozier,~D.~W., Boisvert,~R.~F., Clark,~C.~W.} NIST Handbook of Mathematical Functions. Cambridge University Press, New York, NY, 2010.

\bibitem{Gibbs-phenomenon} \url{https://en.wikipedia.org/wiki/Gibbs_phenomenon}.

\bibitem{Julia_lang} \url{https://en.wikipedia.org/wiki/Julia_(programming_language)},\\
 \url{https://julialang.org/}.


\end{thebibliography}


\appendix

\section{Numerically stable computation of Hahn polynomials}
\label{App-Hahn-computation}


Here we expose our numerically stable method for construction of discrete orthogonal polynomials of higher degrees. It makes possible reliable construction of Hahn-Chebyshev polynomials using standard accuracy (double precision, 8-byte) on thousands of points, for degrees up to several hundred.

The standard methods for construction of Hahn-Chebyshev polynomials are as follows:
\begin{itemize}
  \item by explicit formulas (\ref{Qn}), (\ref{Qn-normalized});
  \item by recurrence relations \cite{NSU85};
  \item by the standard Gram-Schmidt orthogonalization;
  \item by multiple precision calculations. 
\end{itemize}
The first three turn out to be very unstable even for small degree polynomials as one can easily check. 

The multiple precision calculations will produce a reliable answer only if one uses very large precision settings (up to a few hundreds of decimal digits) so the calculation speed will reduce dramatically.

The proposed algorithm 
is capable to build discrete orthogonal polynomials for lattices with a (slightly) uneven time step (missing values of the time series, etc.). One should keep this in mind when analyzing the steps of the algorithm below. The main idea consists in computing the values of Legendre polynomials on the input lattice \(\{x_j\}\) and application of an appropriately modified Gram-Schmidt orthogonalization.


The tested version of Julia programming language \cite{Julia_lang} implementation of our algorithm is available at\\
\url{https://github.com/sptsarev/high-deg-polynomial-fitting}.

The Julia scripts posted there include 4 version of the algorithm: for equidistant/non-equidistant lattice and for computation with standard/extended precision.
The scripts are ready to use for fitting time series data with best least squares high-degree polynomial approximation. 

Below we give the main part of our scripts: construction of the normalized Hahn-Chebyshev univariate polynomials (in a form of some commented metacode).

\begin{algorithm}
  \KwIn{List of the lattice points \(\{x_j\}\), \(j=0,\ldots,N\);
   maximal degree \(M\) of the polynomials to be constructed; 
   orthogonalization tolerance level $orth\_tol$.
  }
  \KwOut{$(N+1)\times(M+1)$-matrix of values \(\Q_{m}^{N}(x_j)\) of the normalized Hahn-Chebyshev univariate polynomials of degrees \(m=0,\ldots,M\) on the input lattice \(\{x_j\}\).}
  \Begin{
    \tcc{Everywhere below $j$ denotes the array index \emph{implicitly} running from 0 to $N$.}
    \tcc{Array $x[j]$ contains the points of the lattice  \(\{x_{j}\}\).}
    $x[j] := $READ(data file with \(\{x_{j}\}\))\;
    \tcc{normalize the interval \([x_0 , x_N]\) to $[-1,1]$:}
    $x0 := x[0]$;\quad 
    $xN := x[N]$;\quad 
    $xn[j] := 2  (x[j] - x0) / (xN - x0) - 1$\;
    \tcc{Matrix element $Q[j,m]$ will contain the value of  \(\Q_{m}^{N}(x_{j})\):}
    $Q = Matrix(j = 0 \ldots N, m = 0 \ldots M)$\;
    \tcc{Computing Legendre polynomials on the grid \(\{\hat x_j\}\in [-1,1]\). \\
    They are NOT orthogonal on this grid!}
    $Q[j,0] := 1.0$;\quad 
    $Q[j,1] := xn[j]$\;
    \For{ $d = 2 \ldots M$}{
        $n := d - 1$\;
        $Q[j,d] :=   ((2n + 1)\cdot xn[j]\cdot Q[j,d - 1] - n\cdot Q[j,d - 2]) / (n + 1)$\;
    }
    \tcc{Orthogonalization and normalization of the constructed Legendre polynomials on the input lattice:}
    $norm := \sqrt{\sum_j Q[j,0]^2}$; \quad    $Q[j,0] := Q[j,0] / norm $\;
    \For{ $n = 1 \ldots M$}{
        $u[j] := Q[j,n]$\;
        $orth\_err := \infty $\;
        \While{ $orth\_err > 2 N \cdot orth\_tol $}{
            \For{ $i = 0 \ldots (n - 1)$}{
                $c := \sum_j(u[j]\cdot Q[j,i])$\;
                $u[j] := u[j] - c \cdot Q[j,i]$\; 
            }
            $orth\_err := 0.0$\;
            \For{ $i = 0 \ldots (n - 1)$}{
                $orth\_err := \max(orth\_err, |\sum_j (u[j]\cdot Q[j,i])|)$\;
            }
        }
        $norm := \sqrt{\sum_j u[j]^2}$\;
        $Q[j,n] := u[j] / norm$\;
    }
  \KwRet{Matrix $ Q[j,m] $}
  }
  \caption{Practical robust algorithm for construction of Hahn-Chebyshev univariate polynomials with unit weights on uniform and slightly non-uniform point lattices.\label{alg_1}}
  \end{algorithm}


\end{document}